\documentclass[12pt]{amsart}

\usepackage{amsmath}
\usepackage{amsfonts}
\usepackage{amssymb}
\usepackage{amsthm}
\usepackage{url}
\usepackage{enumerate}

\usepackage[latin1]{inputenc}

\newtheorem{theorem}{Theorem}
\newtheorem{lemma}[theorem]{Lemma}

\newtheorem{corollary}[theorem]{Corollary}

\theoremstyle{definition}

\newtheorem{remark}[theorem]{Remark}

\newcommand{\Z}{{\mathbb Z}}

\newcommand{\R}{\mathbb R}
\newcommand{\Hilbert}{\mathcal H}

\newcommand{\Normaliz}{{\tt Normaliz}}
\newcommand{\boproof}{\noindent {\bf Proof. }}
\newcommand{\eoproof}{\hspace*{\fill} $\square$ \vspace{5pt}}

\DeclareMathOperator{\cone}{cone}
\DeclareMathOperator{\monoid}{monoid}
\DeclareMathOperator{\lattice}{lattice}
\def\cA{{\mathcal A}}

\usepackage{enumerate}
\def\ve#1{\mathchoice{\mbox{\boldmath$\displaystyle\bf#1$}}
{\mbox{\boldmath$\textstyle\bf#1$}}
{\mbox{\boldmath$\scriptstyle\bf#1$}}
{\mbox{\boldmath$\scriptscriptstyle\bf#1$}}}

\newcommand\vecc{{\ve c}}

\newcommand\vef{{\ve f}}
\newcommand\veg{{\ve g}}
\newcommand\veh{{\ve h}}

\newcommand\ver{{\ve r}}
\newcommand\ves{{\ve s}}

\newcommand\vev{{\ve v}}

\newcommand\vex{{\ve x}}

\newcommand\vez{{\ve z}}

\begin{document}

\title[Challenging computations of Hilbert bases]{Challenging computations of Hilbert bases of cones associated with algebraic statistics}

\author[W. Bruns]{Winfried Bruns}
\address{Winfried Bruns\\ Universit\"at Osnabr\"uck\\ FB Mathematik/Informatik\\ 49069 Osna\-br\"uck\\ Germany}
\email{wbruns@uos.de}
\author[R. Hemmecke]{Raymond Hemmecke}
\address{Raymond Hemmecke\\ Zentrum Mathematik, M9\\ Technische Universität M\"unchen\\ Boltzmannstr. 3\\ 85747 Garching\\ Germany}
\email{hemmecke@ma.tum.de}
\author[B. Ichim]{Bogdan Ichim}
\address{Bogdan Ichim\\ Institute of Mathematics\\ C.P. 1-764\\ 70700 Bucharest\\ Romania}
\email{bogdan\_ichim@yahoo.com}
\author[M. K\"oppe]{Matthias K\"oppe}
\address{Matthias K\"oppe\\ University of California, Davis\\ Department of Mathematics\\ One Shields Avenue\\ Davis, CA 95616\\ USA}
\email{mkoeppe@math.ucdavis.edu}
\author[Ch. S\"oger]{Christof S\"oger}
\address{Christof S\"oger\\ Universit\"at Osnabr\"uck\\ FB Mathematik/Informatik\\ 49069 Osna\-br\"uck\\ Germany}
\email{csoeger@uos.de}

\date{\today}

\begin{abstract}
In this paper we present two independent computational proofs
that the monoid derived from $5\times 5\times 3$ contingency
tables is normal, completing the classification by Hibi and
Ohsugi. We show that Vlach's vector disproving normality for
the monoid derived from $6\times 4\times 3$ contingency tables
is the unique minimal such vector up to symmetry. Finally, we
compute the full Hilbert basis of the cone associated with the
non-normal monoid of the semi-graphoid for $|N|=5$. The
computations are based on extensions of the packages LattE-4ti2
and Normaliz.
\end{abstract}

\maketitle

\section{Introduction}
\label{Section: Introduction}
\let\thefootnote\relax\footnotetext{Acknowledgement: B. Icim was partially supported by CNCSIS grant RP-1 no. 7/01.07.2009
during the preparation of this work.}

Let $S=\monoid(G)$ be an affine monoid generated by a finite
set $G\subseteq\Z^n$ of integer vectors. We call $S$
\emph{normal} if $S=\cone(G)\cap\lattice(G)$, where
$\cone(G)=\{\vex\in\R^n:\vex=\sum\lambda_i
\veg_i,\lambda_i\in\R_+,\veg_i\in G\}$ denotes the rational
polyhedral cone generated by $G$ and where
$\lattice(G)=\{\vex\in\R^n:\vex=\sum\lambda_i
\veg_i,\lambda_i\in\Z,\veg_i\in G\}$ denotes the sublattice of
$\Z^n$ generated by $G$. In this paper, we will stick to the
case that $\lattice(G)=\Z^n$. Then, normality of $S$ is
equivalent to saying that $G$ contains the Hilbert basis of
$\cone(G)$, i.e., every lattice point in $\cone(G)$ can be
written as a nonnegative integer linear combination of elements
in $G$. By the \emph{Hilbert basis} $\Hilbert(C)$ of a pointed
rational cone $C$ we mean the unique minimal system of
generators of the monoid $M$ of lattice points in $C$. The
Hilbert basis of $C$ consists of the \emph{irreducible
elements} of $M$, i.e., those elements of $M$ that do not have
a nontrivial representation as a sum of two elements of $M$
(see \cite[Ch.~2]{BG} for a comprehensive discussion).  Note
that deciding normality of an affine monoid is NP-hard
\cite{Durand+Hermann+Juban:99}.

Normality of monoids derived from  $r_1\times
r_2\times\cdots\times r_N$ contingency tables by taking
$N-1$-marginals (that is, line sums) was settled almost
completely by Hibi and Ohsugi \cite{Hibi+Ohsugi}. In this paper
we close the last open cases by showing computationally, via
two different approaches and independent implementations, that
$5\times 5\times 3$ has a normal monoid. The normality for
$5\times5\times3$ implies normality for the other two open
cases $5\times4\times3$ and $4\times4\times3$ by
\cite[3.2]{Hibi+Ohsugi}.

Here is the defining matrix $A_{5\times 5\times 3}$ whose
columns generate the monoid associated to $5\times 5\times 3$
contingency tables. Every $\cdot$ corresponds to an entry $0$.

{\tiny
\[
\begin{array}{c}
1\cdot\cdot\cdot\cdot\cdot\cdot\cdot\cdot\cdot\cdot\cdot\cdot\cdot\cdot1\cdot\cdot\cdot\cdot\cdot\cdot\cdot\cdot\cdot\cdot\cdot\cdot\cdot\cdot1\cdot\cdot\cdot\cdot\cdot\cdot\cdot\cdot\cdot\cdot\cdot\cdot\cdot\cdot1\cdot\cdot\cdot\cdot\cdot\cdot\cdot\cdot\cdot\cdot\cdot\cdot\cdot\cdot1\cdot\cdot\cdot\cdot\cdot\cdot\cdot\cdot\cdot\cdot\cdot\cdot\cdot\cdot\\
\cdot1\cdot\cdot\cdot\cdot\cdot\cdot\cdot\cdot\cdot\cdot\cdot\cdot\cdot\cdot1\cdot\cdot\cdot\cdot\cdot\cdot\cdot\cdot\cdot\cdot\cdot\cdot\cdot\cdot1\cdot\cdot\cdot\cdot\cdot\cdot\cdot\cdot\cdot\cdot\cdot\cdot\cdot\cdot1\cdot\cdot\cdot\cdot\cdot\cdot\cdot\cdot\cdot\cdot\cdot\cdot\cdot\cdot1\cdot\cdot\cdot\cdot\cdot\cdot\cdot\cdot\cdot\cdot\cdot\cdot\cdot\\
\cdot\cdot1\cdot\cdot\cdot\cdot\cdot\cdot\cdot\cdot\cdot\cdot\cdot\cdot\cdot\cdot1\cdot\cdot\cdot\cdot\cdot\cdot\cdot\cdot\cdot\cdot\cdot\cdot\cdot\cdot1\cdot\cdot\cdot\cdot\cdot\cdot\cdot\cdot\cdot\cdot\cdot\cdot\cdot\cdot1\cdot\cdot\cdot\cdot\cdot\cdot\cdot\cdot\cdot\cdot\cdot\cdot\cdot\cdot1\cdot\cdot\cdot\cdot\cdot\cdot\cdot\cdot\cdot\cdot\cdot\cdot\\
\cdot\cdot\cdot1\cdot\cdot\cdot\cdot\cdot\cdot\cdot\cdot\cdot\cdot\cdot\cdot\cdot\cdot1\cdot\cdot\cdot\cdot\cdot\cdot\cdot\cdot\cdot\cdot\cdot\cdot\cdot\cdot1\cdot\cdot\cdot\cdot\cdot\cdot\cdot\cdot\cdot\cdot\cdot\cdot\cdot\cdot1\cdot\cdot\cdot\cdot\cdot\cdot\cdot\cdot\cdot\cdot\cdot\cdot\cdot\cdot1\cdot\cdot\cdot\cdot\cdot\cdot\cdot\cdot\cdot\cdot\cdot\\
\cdot\cdot\cdot\cdot1\cdot\cdot\cdot\cdot\cdot\cdot\cdot\cdot\cdot\cdot\cdot\cdot\cdot\cdot1\cdot\cdot\cdot\cdot\cdot\cdot\cdot\cdot\cdot\cdot\cdot\cdot\cdot\cdot1\cdot\cdot\cdot\cdot\cdot\cdot\cdot\cdot\cdot\cdot\cdot\cdot\cdot\cdot1\cdot\cdot\cdot\cdot\cdot\cdot\cdot\cdot\cdot\cdot\cdot\cdot\cdot\cdot1\cdot\cdot\cdot\cdot\cdot\cdot\cdot\cdot\cdot\cdot\\
\cdot\cdot\cdot\cdot\cdot1\cdot\cdot\cdot\cdot\cdot\cdot\cdot\cdot\cdot\cdot\cdot\cdot\cdot\cdot1\cdot\cdot\cdot\cdot\cdot\cdot\cdot\cdot\cdot\cdot\cdot\cdot\cdot\cdot1\cdot\cdot\cdot\cdot\cdot\cdot\cdot\cdot\cdot\cdot\cdot\cdot\cdot\cdot1\cdot\cdot\cdot\cdot\cdot\cdot\cdot\cdot\cdot\cdot\cdot\cdot\cdot\cdot1\cdot\cdot\cdot\cdot\cdot\cdot\cdot\cdot\cdot\\
\cdot\cdot\cdot\cdot\cdot\cdot1\cdot\cdot\cdot\cdot\cdot\cdot\cdot\cdot\cdot\cdot\cdot\cdot\cdot\cdot1\cdot\cdot\cdot\cdot\cdot\cdot\cdot\cdot\cdot\cdot\cdot\cdot\cdot\cdot1\cdot\cdot\cdot\cdot\cdot\cdot\cdot\cdot\cdot\cdot\cdot\cdot\cdot\cdot1\cdot\cdot\cdot\cdot\cdot\cdot\cdot\cdot\cdot\cdot\cdot\cdot\cdot\cdot1\cdot\cdot\cdot\cdot\cdot\cdot\cdot\cdot\\
\cdot\cdot\cdot\cdot\cdot\cdot\cdot1\cdot\cdot\cdot\cdot\cdot\cdot\cdot\cdot\cdot\cdot\cdot\cdot\cdot\cdot1\cdot\cdot\cdot\cdot\cdot\cdot\cdot\cdot\cdot\cdot\cdot\cdot\cdot\cdot1\cdot\cdot\cdot\cdot\cdot\cdot\cdot\cdot\cdot\cdot\cdot\cdot\cdot\cdot1\cdot\cdot\cdot\cdot\cdot\cdot\cdot\cdot\cdot\cdot\cdot\cdot\cdot\cdot1\cdot\cdot\cdot\cdot\cdot\cdot\cdot\\
\cdot\cdot\cdot\cdot\cdot\cdot\cdot\cdot1\cdot\cdot\cdot\cdot\cdot\cdot\cdot\cdot\cdot\cdot\cdot\cdot\cdot\cdot1\cdot\cdot\cdot\cdot\cdot\cdot\cdot\cdot\cdot\cdot\cdot\cdot\cdot\cdot1\cdot\cdot\cdot\cdot\cdot\cdot\cdot\cdot\cdot\cdot\cdot\cdot\cdot\cdot1\cdot\cdot\cdot\cdot\cdot\cdot\cdot\cdot\cdot\cdot\cdot\cdot\cdot\cdot1\cdot\cdot\cdot\cdot\cdot\cdot\\
\cdot\cdot\cdot\cdot\cdot\cdot\cdot\cdot\cdot1\cdot\cdot\cdot\cdot\cdot\cdot\cdot\cdot\cdot\cdot\cdot\cdot\cdot\cdot1\cdot\cdot\cdot\cdot\cdot\cdot\cdot\cdot\cdot\cdot\cdot\cdot\cdot\cdot1\cdot\cdot\cdot\cdot\cdot\cdot\cdot\cdot\cdot\cdot\cdot\cdot\cdot\cdot1\cdot\cdot\cdot\cdot\cdot\cdot\cdot\cdot\cdot\cdot\cdot\cdot\cdot\cdot1\cdot\cdot\cdot\cdot\cdot\\
\cdot\cdot\cdot\cdot\cdot\cdot\cdot\cdot\cdot\cdot1\cdot\cdot\cdot\cdot\cdot\cdot\cdot\cdot\cdot\cdot\cdot\cdot\cdot\cdot1\cdot\cdot\cdot\cdot\cdot\cdot\cdot\cdot\cdot\cdot\cdot\cdot\cdot\cdot1\cdot\cdot\cdot\cdot\cdot\cdot\cdot\cdot\cdot\cdot\cdot\cdot\cdot\cdot1\cdot\cdot\cdot\cdot\cdot\cdot\cdot\cdot\cdot\cdot\cdot\cdot\cdot\cdot1\cdot\cdot\cdot\cdot\\
\cdot\cdot\cdot\cdot\cdot\cdot\cdot\cdot\cdot\cdot\cdot1\cdot\cdot\cdot\cdot\cdot\cdot\cdot\cdot\cdot\cdot\cdot\cdot\cdot\cdot1\cdot\cdot\cdot\cdot\cdot\cdot\cdot\cdot\cdot\cdot\cdot\cdot\cdot\cdot1\cdot\cdot\cdot\cdot\cdot\cdot\cdot\cdot\cdot\cdot\cdot\cdot\cdot\cdot1\cdot\cdot\cdot\cdot\cdot\cdot\cdot\cdot\cdot\cdot\cdot\cdot\cdot\cdot1\cdot\cdot\cdot\\
\cdot\cdot\cdot\cdot\cdot\cdot\cdot\cdot\cdot\cdot\cdot\cdot1\cdot\cdot\cdot\cdot\cdot\cdot\cdot\cdot\cdot\cdot\cdot\cdot\cdot\cdot1\cdot\cdot\cdot\cdot\cdot\cdot\cdot\cdot\cdot\cdot\cdot\cdot\cdot\cdot1\cdot\cdot\cdot\cdot\cdot\cdot\cdot\cdot\cdot\cdot\cdot\cdot\cdot\cdot1\cdot\cdot\cdot\cdot\cdot\cdot\cdot\cdot\cdot\cdot\cdot\cdot\cdot\cdot1\cdot\cdot\\
\cdot\cdot\cdot\cdot\cdot\cdot\cdot\cdot\cdot\cdot\cdot\cdot\cdot1\cdot\cdot\cdot\cdot\cdot\cdot\cdot\cdot\cdot\cdot\cdot\cdot\cdot\cdot1\cdot\cdot\cdot\cdot\cdot\cdot\cdot\cdot\cdot\cdot\cdot\cdot\cdot\cdot1\cdot\cdot\cdot\cdot\cdot\cdot\cdot\cdot\cdot\cdot\cdot\cdot\cdot\cdot1\cdot\cdot\cdot\cdot\cdot\cdot\cdot\cdot\cdot\cdot\cdot\cdot\cdot\cdot1\cdot\\
\cdot\cdot\cdot\cdot\cdot\cdot\cdot\cdot\cdot\cdot\cdot\cdot\cdot\cdot1\cdot\cdot\cdot\cdot\cdot\cdot\cdot\cdot\cdot\cdot\cdot\cdot\cdot\cdot1\cdot\cdot\cdot\cdot\cdot\cdot\cdot\cdot\cdot\cdot\cdot\cdot\cdot\cdot1\cdot\cdot\cdot\cdot\cdot\cdot\cdot\cdot\cdot\cdot\cdot\cdot\cdot\cdot1\cdot\cdot\cdot\cdot\cdot\cdot\cdot\cdot\cdot\cdot\cdot\cdot\cdot\cdot1\\
1\cdot\cdot1\cdot\cdot1\cdot\cdot1\cdot\cdot1\cdot\cdot\cdot\cdot\cdot\cdot\cdot\cdot\cdot\cdot\cdot\cdot\cdot\cdot\cdot\cdot\cdot\cdot\cdot\cdot\cdot\cdot\cdot\cdot\cdot\cdot\cdot\cdot\cdot\cdot\cdot\cdot\cdot\cdot\cdot\cdot\cdot\cdot\cdot\cdot\cdot\cdot\cdot\cdot\cdot\cdot\cdot\cdot\cdot\cdot\cdot\cdot\cdot\cdot\cdot\cdot\cdot\cdot\cdot\cdot\cdot\cdot\\
\cdot1\cdot\cdot1\cdot\cdot1\cdot\cdot1\cdot\cdot1\cdot\cdot\cdot\cdot\cdot\cdot\cdot\cdot\cdot\cdot\cdot\cdot\cdot\cdot\cdot\cdot\cdot\cdot\cdot\cdot\cdot\cdot\cdot\cdot\cdot\cdot\cdot\cdot\cdot\cdot\cdot\cdot\cdot\cdot\cdot\cdot\cdot\cdot\cdot\cdot\cdot\cdot\cdot\cdot\cdot\cdot\cdot\cdot\cdot\cdot\cdot\cdot\cdot\cdot\cdot\cdot\cdot\cdot\cdot\cdot\cdot\\
\cdot\cdot1\cdot\cdot1\cdot\cdot1\cdot\cdot1\cdot\cdot1\cdot\cdot\cdot\cdot\cdot\cdot\cdot\cdot\cdot\cdot\cdot\cdot\cdot\cdot\cdot\cdot\cdot\cdot\cdot\cdot\cdot\cdot\cdot\cdot\cdot\cdot\cdot\cdot\cdot\cdot\cdot\cdot\cdot\cdot\cdot\cdot\cdot\cdot\cdot\cdot\cdot\cdot\cdot\cdot\cdot\cdot\cdot\cdot\cdot\cdot\cdot\cdot\cdot\cdot\cdot\cdot\cdot\cdot\cdot\cdot\\
\cdot\cdot\cdot\cdot\cdot\cdot\cdot\cdot\cdot\cdot\cdot\cdot\cdot\cdot\cdot1\cdot\cdot1\cdot\cdot1\cdot\cdot1\cdot\cdot1\cdot\cdot\cdot\cdot\cdot\cdot\cdot\cdot\cdot\cdot\cdot\cdot\cdot\cdot\cdot\cdot\cdot\cdot\cdot\cdot\cdot\cdot\cdot\cdot\cdot\cdot\cdot\cdot\cdot\cdot\cdot\cdot\cdot\cdot\cdot\cdot\cdot\cdot\cdot\cdot\cdot\cdot\cdot\cdot\cdot\cdot\cdot\\
\cdot\cdot\cdot\cdot\cdot\cdot\cdot\cdot\cdot\cdot\cdot\cdot\cdot\cdot\cdot\cdot1\cdot\cdot1\cdot\cdot1\cdot\cdot1\cdot\cdot1\cdot\cdot\cdot\cdot\cdot\cdot\cdot\cdot\cdot\cdot\cdot\cdot\cdot\cdot\cdot\cdot\cdot\cdot\cdot\cdot\cdot\cdot\cdot\cdot\cdot\cdot\cdot\cdot\cdot\cdot\cdot\cdot\cdot\cdot\cdot\cdot\cdot\cdot\cdot\cdot\cdot\cdot\cdot\cdot\cdot\cdot\\
\cdot\cdot\cdot\cdot\cdot\cdot\cdot\cdot\cdot\cdot\cdot\cdot\cdot\cdot\cdot\cdot\cdot1\cdot\cdot1\cdot\cdot1\cdot\cdot1\cdot\cdot1\cdot\cdot\cdot\cdot\cdot\cdot\cdot\cdot\cdot\cdot\cdot\cdot\cdot\cdot\cdot\cdot\cdot\cdot\cdot\cdot\cdot\cdot\cdot\cdot\cdot\cdot\cdot\cdot\cdot\cdot\cdot\cdot\cdot\cdot\cdot\cdot\cdot\cdot\cdot\cdot\cdot\cdot\cdot\cdot\cdot\\
\cdot\cdot\cdot\cdot\cdot\cdot\cdot\cdot\cdot\cdot\cdot\cdot\cdot\cdot\cdot\cdot\cdot\cdot\cdot\cdot\cdot\cdot\cdot\cdot\cdot\cdot\cdot\cdot\cdot\cdot1\cdot\cdot1\cdot\cdot1\cdot\cdot1\cdot\cdot1\cdot\cdot\cdot\cdot\cdot\cdot\cdot\cdot\cdot\cdot\cdot\cdot\cdot\cdot\cdot\cdot\cdot\cdot\cdot\cdot\cdot\cdot\cdot\cdot\cdot\cdot\cdot\cdot\cdot\cdot\cdot\cdot\\
\cdot\cdot\cdot\cdot\cdot\cdot\cdot\cdot\cdot\cdot\cdot\cdot\cdot\cdot\cdot\cdot\cdot\cdot\cdot\cdot\cdot\cdot\cdot\cdot\cdot\cdot\cdot\cdot\cdot\cdot\cdot1\cdot\cdot1\cdot\cdot1\cdot\cdot1\cdot\cdot1\cdot\cdot\cdot\cdot\cdot\cdot\cdot\cdot\cdot\cdot\cdot\cdot\cdot\cdot\cdot\cdot\cdot\cdot\cdot\cdot\cdot\cdot\cdot\cdot\cdot\cdot\cdot\cdot\cdot\cdot\cdot\\
\cdot\cdot\cdot\cdot\cdot\cdot\cdot\cdot\cdot\cdot\cdot\cdot\cdot\cdot\cdot\cdot\cdot\cdot\cdot\cdot\cdot\cdot\cdot\cdot\cdot\cdot\cdot\cdot\cdot\cdot\cdot\cdot1\cdot\cdot1\cdot\cdot1\cdot\cdot1\cdot\cdot1\cdot\cdot\cdot\cdot\cdot\cdot\cdot\cdot\cdot\cdot\cdot\cdot\cdot\cdot\cdot\cdot\cdot\cdot\cdot\cdot\cdot\cdot\cdot\cdot\cdot\cdot\cdot\cdot\cdot\cdot\\
\cdot\cdot\cdot\cdot\cdot\cdot\cdot\cdot\cdot\cdot\cdot\cdot\cdot\cdot\cdot\cdot\cdot\cdot\cdot\cdot\cdot\cdot\cdot\cdot\cdot\cdot\cdot\cdot\cdot\cdot\cdot\cdot\cdot\cdot\cdot\cdot\cdot\cdot\cdot\cdot\cdot\cdot\cdot\cdot\cdot1\cdot\cdot1\cdot\cdot1\cdot\cdot1\cdot\cdot1\cdot\cdot\cdot\cdot\cdot\cdot\cdot\cdot\cdot\cdot\cdot\cdot\cdot\cdot\cdot\cdot\cdot\\
\cdot\cdot\cdot\cdot\cdot\cdot\cdot\cdot\cdot\cdot\cdot\cdot\cdot\cdot\cdot\cdot\cdot\cdot\cdot\cdot\cdot\cdot\cdot\cdot\cdot\cdot\cdot\cdot\cdot\cdot\cdot\cdot\cdot\cdot\cdot\cdot\cdot\cdot\cdot\cdot\cdot\cdot\cdot\cdot\cdot\cdot1\cdot\cdot1\cdot\cdot1\cdot\cdot1\cdot\cdot1\cdot\cdot\cdot\cdot\cdot\cdot\cdot\cdot\cdot\cdot\cdot\cdot\cdot\cdot\cdot\cdot\\
\cdot\cdot\cdot\cdot\cdot\cdot\cdot\cdot\cdot\cdot\cdot\cdot\cdot\cdot\cdot\cdot\cdot\cdot\cdot\cdot\cdot\cdot\cdot\cdot\cdot\cdot\cdot\cdot\cdot\cdot\cdot\cdot\cdot\cdot\cdot\cdot\cdot\cdot\cdot\cdot\cdot\cdot\cdot\cdot\cdot\cdot\cdot1\cdot\cdot1\cdot\cdot1\cdot\cdot1\cdot\cdot1\cdot\cdot\cdot\cdot\cdot\cdot\cdot\cdot\cdot\cdot\cdot\cdot\cdot\cdot\cdot\\
\cdot\cdot\cdot\cdot\cdot\cdot\cdot\cdot\cdot\cdot\cdot\cdot\cdot\cdot\cdot\cdot\cdot\cdot\cdot\cdot\cdot\cdot\cdot\cdot\cdot\cdot\cdot\cdot\cdot\cdot\cdot\cdot\cdot\cdot\cdot\cdot\cdot\cdot\cdot\cdot\cdot\cdot\cdot\cdot\cdot\cdot\cdot\cdot\cdot\cdot\cdot\cdot\cdot\cdot\cdot\cdot\cdot\cdot\cdot\cdot1\cdot\cdot1\cdot\cdot1\cdot\cdot1\cdot\cdot1\cdot\cdot\\
\cdot\cdot\cdot\cdot\cdot\cdot\cdot\cdot\cdot\cdot\cdot\cdot\cdot\cdot\cdot\cdot\cdot\cdot\cdot\cdot\cdot\cdot\cdot\cdot\cdot\cdot\cdot\cdot\cdot\cdot\cdot\cdot\cdot\cdot\cdot\cdot\cdot\cdot\cdot\cdot\cdot\cdot\cdot\cdot\cdot\cdot\cdot\cdot\cdot\cdot\cdot\cdot\cdot\cdot\cdot\cdot\cdot\cdot\cdot\cdot\cdot1\cdot\cdot1\cdot\cdot1\cdot\cdot1\cdot\cdot1\cdot\\
\cdot\cdot\cdot\cdot\cdot\cdot\cdot\cdot\cdot\cdot\cdot\cdot\cdot\cdot\cdot\cdot\cdot\cdot\cdot\cdot\cdot\cdot\cdot\cdot\cdot\cdot\cdot\cdot\cdot\cdot\cdot\cdot\cdot\cdot\cdot\cdot\cdot\cdot\cdot\cdot\cdot\cdot\cdot\cdot\cdot\cdot\cdot\cdot\cdot\cdot\cdot\cdot\cdot\cdot\cdot\cdot\cdot\cdot\cdot\cdot\cdot\cdot1\cdot\cdot1\cdot\cdot1\cdot\cdot1\cdot\cdot1\\
111\cdot\cdot\cdot\cdot\cdot\cdot\cdot\cdot\cdot\cdot\cdot\cdot\cdot\cdot\cdot\cdot\cdot\cdot\cdot\cdot\cdot\cdot\cdot\cdot\cdot\cdot\cdot\cdot\cdot\cdot\cdot\cdot\cdot\cdot\cdot\cdot\cdot\cdot\cdot\cdot\cdot\cdot\cdot\cdot\cdot\cdot\cdot\cdot\cdot\cdot\cdot\cdot\cdot\cdot\cdot\cdot\cdot\cdot\cdot\cdot\cdot\cdot\cdot\cdot\cdot\cdot\cdot\cdot\cdot\cdot\cdot\cdot\\
\cdot\cdot\cdot111\cdot\cdot\cdot\cdot\cdot\cdot\cdot\cdot\cdot\cdot\cdot\cdot\cdot\cdot\cdot\cdot\cdot\cdot\cdot\cdot\cdot\cdot\cdot\cdot\cdot\cdot\cdot\cdot\cdot\cdot\cdot\cdot\cdot\cdot\cdot\cdot\cdot\cdot\cdot\cdot\cdot\cdot\cdot\cdot\cdot\cdot\cdot\cdot\cdot\cdot\cdot\cdot\cdot\cdot\cdot\cdot\cdot\cdot\cdot\cdot\cdot\cdot\cdot\cdot\cdot\cdot\cdot\cdot\cdot\\
\cdot\cdot\cdot\cdot\cdot\cdot111\cdot\cdot\cdot\cdot\cdot\cdot\cdot\cdot\cdot\cdot\cdot\cdot\cdot\cdot\cdot\cdot\cdot\cdot\cdot\cdot\cdot\cdot\cdot\cdot\cdot\cdot\cdot\cdot\cdot\cdot\cdot\cdot\cdot\cdot\cdot\cdot\cdot\cdot\cdot\cdot\cdot\cdot\cdot\cdot\cdot\cdot\cdot\cdot\cdot\cdot\cdot\cdot\cdot\cdot\cdot\cdot\cdot\cdot\cdot\cdot\cdot\cdot\cdot\cdot\cdot\cdot\\
\cdot\cdot\cdot\cdot\cdot\cdot\cdot\cdot\cdot111\cdot\cdot\cdot\cdot\cdot\cdot\cdot\cdot\cdot\cdot\cdot\cdot\cdot\cdot\cdot\cdot\cdot\cdot\cdot\cdot\cdot\cdot\cdot\cdot\cdot\cdot\cdot\cdot\cdot\cdot\cdot\cdot\cdot\cdot\cdot\cdot\cdot\cdot\cdot\cdot\cdot\cdot\cdot\cdot\cdot\cdot\cdot\cdot\cdot\cdot\cdot\cdot\cdot\cdot\cdot\cdot\cdot\cdot\cdot\cdot\cdot\cdot\cdot\\
\cdot\cdot\cdot\cdot\cdot\cdot\cdot\cdot\cdot\cdot\cdot\cdot111\cdot\cdot\cdot\cdot\cdot\cdot\cdot\cdot\cdot\cdot\cdot\cdot\cdot\cdot\cdot\cdot\cdot\cdot\cdot\cdot\cdot\cdot\cdot\cdot\cdot\cdot\cdot\cdot\cdot\cdot\cdot\cdot\cdot\cdot\cdot\cdot\cdot\cdot\cdot\cdot\cdot\cdot\cdot\cdot\cdot\cdot\cdot\cdot\cdot\cdot\cdot\cdot\cdot\cdot\cdot\cdot\cdot\cdot\cdot\cdot\\
\cdot\cdot\cdot\cdot\cdot\cdot\cdot\cdot\cdot\cdot\cdot\cdot\cdot\cdot\cdot111\cdot\cdot\cdot\cdot\cdot\cdot\cdot\cdot\cdot\cdot\cdot\cdot\cdot\cdot\cdot\cdot\cdot\cdot\cdot\cdot\cdot\cdot\cdot\cdot\cdot\cdot\cdot\cdot\cdot\cdot\cdot\cdot\cdot\cdot\cdot\cdot\cdot\cdot\cdot\cdot\cdot\cdot\cdot\cdot\cdot\cdot\cdot\cdot\cdot\cdot\cdot\cdot\cdot\cdot\cdot\cdot\cdot\\
\cdot\cdot\cdot\cdot\cdot\cdot\cdot\cdot\cdot\cdot\cdot\cdot\cdot\cdot\cdot\cdot\cdot\cdot111\cdot\cdot\cdot\cdot\cdot\cdot\cdot\cdot\cdot\cdot\cdot\cdot\cdot\cdot\cdot\cdot\cdot\cdot\cdot\cdot\cdot\cdot\cdot\cdot\cdot\cdot\cdot\cdot\cdot\cdot\cdot\cdot\cdot\cdot\cdot\cdot\cdot\cdot\cdot\cdot\cdot\cdot\cdot\cdot\cdot\cdot\cdot\cdot\cdot\cdot\cdot\cdot\cdot\cdot\\
\cdot\cdot\cdot\cdot\cdot\cdot\cdot\cdot\cdot\cdot\cdot\cdot\cdot\cdot\cdot\cdot\cdot\cdot\cdot\cdot\cdot111\cdot\cdot\cdot\cdot\cdot\cdot\cdot\cdot\cdot\cdot\cdot\cdot\cdot\cdot\cdot\cdot\cdot\cdot\cdot\cdot\cdot\cdot\cdot\cdot\cdot\cdot\cdot\cdot\cdot\cdot\cdot\cdot\cdot\cdot\cdot\cdot\cdot\cdot\cdot\cdot\cdot\cdot\cdot\cdot\cdot\cdot\cdot\cdot\cdot\cdot\cdot\\
\cdot\cdot\cdot\cdot\cdot\cdot\cdot\cdot\cdot\cdot\cdot\cdot\cdot\cdot\cdot\cdot\cdot\cdot\cdot\cdot\cdot\cdot\cdot\cdot111\cdot\cdot\cdot\cdot\cdot\cdot\cdot\cdot\cdot\cdot\cdot\cdot\cdot\cdot\cdot\cdot\cdot\cdot\cdot\cdot\cdot\cdot\cdot\cdot\cdot\cdot\cdot\cdot\cdot\cdot\cdot\cdot\cdot\cdot\cdot\cdot\cdot\cdot\cdot\cdot\cdot\cdot\cdot\cdot\cdot\cdot\cdot\cdot\\
\cdot\cdot\cdot\cdot\cdot\cdot\cdot\cdot\cdot\cdot\cdot\cdot\cdot\cdot\cdot\cdot\cdot\cdot\cdot\cdot\cdot\cdot\cdot\cdot\cdot\cdot\cdot111\cdot\cdot\cdot\cdot\cdot\cdot\cdot\cdot\cdot\cdot\cdot\cdot\cdot\cdot\cdot\cdot\cdot\cdot\cdot\cdot\cdot\cdot\cdot\cdot\cdot\cdot\cdot\cdot\cdot\cdot\cdot\cdot\cdot\cdot\cdot\cdot\cdot\cdot\cdot\cdot\cdot\cdot\cdot\cdot\cdot\\
\cdot\cdot\cdot\cdot\cdot\cdot\cdot\cdot\cdot\cdot\cdot\cdot\cdot\cdot\cdot\cdot\cdot\cdot\cdot\cdot\cdot\cdot\cdot\cdot\cdot\cdot\cdot\cdot\cdot\cdot111\cdot\cdot\cdot\cdot\cdot\cdot\cdot\cdot\cdot\cdot\cdot\cdot\cdot\cdot\cdot\cdot\cdot\cdot\cdot\cdot\cdot\cdot\cdot\cdot\cdot\cdot\cdot\cdot\cdot\cdot\cdot\cdot\cdot\cdot\cdot\cdot\cdot\cdot\cdot\cdot\cdot\cdot\\
\cdot\cdot\cdot\cdot\cdot\cdot\cdot\cdot\cdot\cdot\cdot\cdot\cdot\cdot\cdot\cdot\cdot\cdot\cdot\cdot\cdot\cdot\cdot\cdot\cdot\cdot\cdot\cdot\cdot\cdot\cdot\cdot\cdot111\cdot\cdot\cdot\cdot\cdot\cdot\cdot\cdot\cdot\cdot\cdot\cdot\cdot\cdot\cdot\cdot\cdot\cdot\cdot\cdot\cdot\cdot\cdot\cdot\cdot\cdot\cdot\cdot\cdot\cdot\cdot\cdot\cdot\cdot\cdot\cdot\cdot\cdot\cdot\\
\cdot\cdot\cdot\cdot\cdot\cdot\cdot\cdot\cdot\cdot\cdot\cdot\cdot\cdot\cdot\cdot\cdot\cdot\cdot\cdot\cdot\cdot\cdot\cdot\cdot\cdot\cdot\cdot\cdot\cdot\cdot\cdot\cdot\cdot\cdot\cdot111\cdot\cdot\cdot\cdot\cdot\cdot\cdot\cdot\cdot\cdot\cdot\cdot\cdot\cdot\cdot\cdot\cdot\cdot\cdot\cdot\cdot\cdot\cdot\cdot\cdot\cdot\cdot\cdot\cdot\cdot\cdot\cdot\cdot\cdot\cdot\cdot\\
\cdot\cdot\cdot\cdot\cdot\cdot\cdot\cdot\cdot\cdot\cdot\cdot\cdot\cdot\cdot\cdot\cdot\cdot\cdot\cdot\cdot\cdot\cdot\cdot\cdot\cdot\cdot\cdot\cdot\cdot\cdot\cdot\cdot\cdot\cdot\cdot\cdot\cdot\cdot111\cdot\cdot\cdot\cdot\cdot\cdot\cdot\cdot\cdot\cdot\cdot\cdot\cdot\cdot\cdot\cdot\cdot\cdot\cdot\cdot\cdot\cdot\cdot\cdot\cdot\cdot\cdot\cdot\cdot\cdot\cdot\cdot\cdot\\
\cdot\cdot\cdot\cdot\cdot\cdot\cdot\cdot\cdot\cdot\cdot\cdot\cdot\cdot\cdot\cdot\cdot\cdot\cdot\cdot\cdot\cdot\cdot\cdot\cdot\cdot\cdot\cdot\cdot\cdot\cdot\cdot\cdot\cdot\cdot\cdot\cdot\cdot\cdot\cdot\cdot\cdot111\cdot\cdot\cdot\cdot\cdot\cdot\cdot\cdot\cdot\cdot\cdot\cdot\cdot\cdot\cdot\cdot\cdot\cdot\cdot\cdot\cdot\cdot\cdot\cdot\cdot\cdot\cdot\cdot\cdot\cdot\\
\cdot\cdot\cdot\cdot\cdot\cdot\cdot\cdot\cdot\cdot\cdot\cdot\cdot\cdot\cdot\cdot\cdot\cdot\cdot\cdot\cdot\cdot\cdot\cdot\cdot\cdot\cdot\cdot\cdot\cdot\cdot\cdot\cdot\cdot\cdot\cdot\cdot\cdot\cdot\cdot\cdot\cdot\cdot\cdot\cdot111\cdot\cdot\cdot\cdot\cdot\cdot\cdot\cdot\cdot\cdot\cdot\cdot\cdot\cdot\cdot\cdot\cdot\cdot\cdot\cdot\cdot\cdot\cdot\cdot\cdot\cdot\cdot\\
\cdot\cdot\cdot\cdot\cdot\cdot\cdot\cdot\cdot\cdot\cdot\cdot\cdot\cdot\cdot\cdot\cdot\cdot\cdot\cdot\cdot\cdot\cdot\cdot\cdot\cdot\cdot\cdot\cdot\cdot\cdot\cdot\cdot\cdot\cdot\cdot\cdot\cdot\cdot\cdot\cdot\cdot\cdot\cdot\cdot\cdot\cdot\cdot111\cdot\cdot\cdot\cdot\cdot\cdot\cdot\cdot\cdot\cdot\cdot\cdot\cdot\cdot\cdot\cdot\cdot\cdot\cdot\cdot\cdot\cdot\cdot\cdot\\
\cdot\cdot\cdot\cdot\cdot\cdot\cdot\cdot\cdot\cdot\cdot\cdot\cdot\cdot\cdot\cdot\cdot\cdot\cdot\cdot\cdot\cdot\cdot\cdot\cdot\cdot\cdot\cdot\cdot\cdot\cdot\cdot\cdot\cdot\cdot\cdot\cdot\cdot\cdot\cdot\cdot\cdot\cdot\cdot\cdot\cdot\cdot\cdot\cdot\cdot\cdot111\cdot\cdot\cdot\cdot\cdot\cdot\cdot\cdot\cdot\cdot\cdot\cdot\cdot\cdot\cdot\cdot\cdot\cdot\cdot\cdot\cdot\\
\cdot\cdot\cdot\cdot\cdot\cdot\cdot\cdot\cdot\cdot\cdot\cdot\cdot\cdot\cdot\cdot\cdot\cdot\cdot\cdot\cdot\cdot\cdot\cdot\cdot\cdot\cdot\cdot\cdot\cdot\cdot\cdot\cdot\cdot\cdot\cdot\cdot\cdot\cdot\cdot\cdot\cdot\cdot\cdot\cdot\cdot\cdot\cdot\cdot\cdot\cdot\cdot\cdot\cdot111\cdot\cdot\cdot\cdot\cdot\cdot\cdot\cdot\cdot\cdot\cdot\cdot\cdot\cdot\cdot\cdot\cdot\cdot\\
\cdot\cdot\cdot\cdot\cdot\cdot\cdot\cdot\cdot\cdot\cdot\cdot\cdot\cdot\cdot\cdot\cdot\cdot\cdot\cdot\cdot\cdot\cdot\cdot\cdot\cdot\cdot\cdot\cdot\cdot\cdot\cdot\cdot\cdot\cdot\cdot\cdot\cdot\cdot\cdot\cdot\cdot\cdot\cdot\cdot\cdot\cdot\cdot\cdot\cdot\cdot\cdot\cdot\cdot\cdot\cdot\cdot111\cdot\cdot\cdot\cdot\cdot\cdot\cdot\cdot\cdot\cdot\cdot\cdot\cdot\cdot\cdot\\
\cdot\cdot\cdot\cdot\cdot\cdot\cdot\cdot\cdot\cdot\cdot\cdot\cdot\cdot\cdot\cdot\cdot\cdot\cdot\cdot\cdot\cdot\cdot\cdot\cdot\cdot\cdot\cdot\cdot\cdot\cdot\cdot\cdot\cdot\cdot\cdot\cdot\cdot\cdot\cdot\cdot\cdot\cdot\cdot\cdot\cdot\cdot\cdot\cdot\cdot\cdot\cdot\cdot\cdot\cdot\cdot\cdot\cdot\cdot\cdot111\cdot\cdot\cdot\cdot\cdot\cdot\cdot\cdot\cdot\cdot\cdot\cdot\\
\cdot\cdot\cdot\cdot\cdot\cdot\cdot\cdot\cdot\cdot\cdot\cdot\cdot\cdot\cdot\cdot\cdot\cdot\cdot\cdot\cdot\cdot\cdot\cdot\cdot\cdot\cdot\cdot\cdot\cdot\cdot\cdot\cdot\cdot\cdot\cdot\cdot\cdot\cdot\cdot\cdot\cdot\cdot\cdot\cdot\cdot\cdot\cdot\cdot\cdot\cdot\cdot\cdot\cdot\cdot\cdot\cdot\cdot\cdot\cdot\cdot\cdot\cdot111\cdot\cdot\cdot\cdot\cdot\cdot\cdot\cdot\cdot\\
\cdot\cdot\cdot\cdot\cdot\cdot\cdot\cdot\cdot\cdot\cdot\cdot\cdot\cdot\cdot\cdot\cdot\cdot\cdot\cdot\cdot\cdot\cdot\cdot\cdot\cdot\cdot\cdot\cdot\cdot\cdot\cdot\cdot\cdot\cdot\cdot\cdot\cdot\cdot\cdot\cdot\cdot\cdot\cdot\cdot\cdot\cdot\cdot\cdot\cdot\cdot\cdot\cdot\cdot\cdot\cdot\cdot\cdot\cdot\cdot\cdot\cdot\cdot\cdot\cdot\cdot111\cdot\cdot\cdot\cdot\cdot\cdot\\
\cdot\cdot\cdot\cdot\cdot\cdot\cdot\cdot\cdot\cdot\cdot\cdot\cdot\cdot\cdot\cdot\cdot\cdot\cdot\cdot\cdot\cdot\cdot\cdot\cdot\cdot\cdot\cdot\cdot\cdot\cdot\cdot\cdot\cdot\cdot\cdot\cdot\cdot\cdot\cdot\cdot\cdot\cdot\cdot\cdot\cdot\cdot\cdot\cdot\cdot\cdot\cdot\cdot\cdot\cdot\cdot\cdot\cdot\cdot\cdot\cdot\cdot\cdot\cdot\cdot\cdot\cdot\cdot\cdot111\cdot\cdot\cdot\\
\cdot\cdot\cdot\cdot\cdot\cdot\cdot\cdot\cdot\cdot\cdot\cdot\cdot\cdot\cdot\cdot\cdot\cdot\cdot\cdot\cdot\cdot\cdot\cdot\cdot\cdot\cdot\cdot\cdot\cdot\cdot\cdot\cdot\cdot\cdot\cdot\cdot\cdot\cdot\cdot\cdot\cdot\cdot\cdot\cdot\cdot\cdot\cdot\cdot\cdot\cdot\cdot\cdot\cdot\cdot\cdot\cdot\cdot\cdot\cdot\cdot\cdot\cdot\cdot\cdot\cdot\cdot\cdot\cdot\cdot\cdot\cdot111\\
\end{array}
\]
}

Note that this normality problem cannot be settled directly by
computing the Hilbert basis of the associated cone using
state-of-the-art software such as {\tt Normaliz} v2.2
\cite{Normaliz, BI} or {\tt 4ti2} v1.3.2 \cite{4ti2,
Hemmecke:Hilbert}. Both codes fail to return an answer due to
time and to memory requirements of intermediate computations.
Using the computational approaches presented below, we can now
show the following.

\begin{lemma}\label{Lemma: 355 is normal}
The monoid derived from of $5\times 5\times 3$ contingency
tables by taking line sums (= two-marginals) is normal.
\end{lemma}

This completes the normality classification of the monoids
derived from $r_1\times r_2\times\cdots\times r_N$ contingency
tables by taking line sums as given in \cite{Hibi+Ohsugi}:

\begin{theorem}
Let $r_1\geq r_2\geq\ldots\geq r_N\geq 2$ be integer numbers.
Then the monoid derived from $r_1\times r_2\times\cdots\times
r_N$ contingency tables by taking line sums is normal if and
only if the contingency table is of size
\begin{itemize}
  \item $r_1\times r_2$, $r_1\times r_2\times
      2\times\ldots\times 2$, or
  \item $r_1\times 3\times 3$, or
  \item $4\times 4\times 3$, $5\times 4\times 3$, or
      $5\times 5\times 3$.
\end{itemize}
\end{theorem}

For the monoid of $6\times 4\times 3$ contingency tables, a
vector disproving normality was presented by Vlach
\cite{Vlach:86}. The right-hand side vector $\vef$ for the
counts along the coordinate axes is given by the following
three matrices:
\[
\begin{pmatrix}
1 & 1 & 1\\
1 & 1 & 1\\
1 & 1 & 1\\
1 & 1 & 1\\
\end{pmatrix},
\begin{pmatrix}
1 & 1 & 0\\
1 & 1 & 0\\
1 & 0 & 1\\
1 & 0 & 1\\
0 & 1 & 1\\
0 & 1 & 1\\
\end{pmatrix} \text{ and }
\begin{pmatrix}
1 & 0 & 0 & 1\\
0 & 1 & 1 & 0\\
1 & 1 & 0 & 0\\
0 & 0 & 1 & 1\\
1 & 0 & 1 & 0\\
0 & 1 & 0 & 1\\
\end{pmatrix}.
\]
The unique point in the $6\times 4\times 3$ transportation
polytope $\{\vez\in\R^{72}:A\vez=\vef, \vez\geq\ve 0\}$ is
\[
\vez^*=\frac{1}{2}\left(
\begin{array}{c|c|c|c|c|c}
\begin{array}{ccc}
1 & 1 & 0 \\
0 & 0 & 0 \\
0 & 0 & 0 \\
1 & 1 & 0 \\
\end{array} &
\begin{array}{ccc}
0 & 0 & 0 \\
1 & 1 & 0 \\
1 & 1 & 0 \\
0 & 0 & 0 \\
\end{array} &
\begin{array}{ccc}
1 & 0 & 1 \\
1 & 0 & 1 \\
0 & 0 & 0 \\
0 & 0 & 0 \\
\end{array} &
\begin{array}{ccc}
0 & 0 & 0 \\
0 & 0 & 0 \\
1 & 0 & 1 \\
1 & 0 & 1 \\
\end{array} &
\begin{array}{ccc}
0 & 1 & 1 \\
0 & 0 & 0 \\
0 & 1 & 1 \\
0 & 0 & 0 \\
\end{array} &
\begin{array}{ccc}
0 & 0 & 0 \\
0 & 1 & 1 \\
0 & 0 & 0 \\
0 & 1 & 1\\
\end{array}
\end{array}
\right).
\]
So $\vez^*$ is indeed a hole of the $6\times 4\times 3$ monoid.
We are able to show the following.

\begin{lemma}
The right-hand side vector $\vef$ presented by Vlach
\cite{Vlach:86} is the unique vector (up to the underlying
$S_6\times S_4\times S_3$ symmetry) in the Hilbert basis of the
cone of $6\times 4\times 3$ contingency tables that is not an
extreme ray.
\end{lemma}

The treatment in \cite{Hemmecke+Takemura+Yoshida} now
completely describes \emph{all} holes of the cone, that is, all
lattice points in $\cone(A_{6\times 4\times 3})$ that cannot be
written as a nonnegative linear integer combination of the
(integer) generators of the cone:

\begin{corollary}
Let $\vef$ be the hole in $\cone(A_{6\times 4\times 3})$
presented by Vlach \cite{Vlach:86} and let $\vez^*\in\R^{72}_+$
be the unique solution to $A_{6\times 4\times 3} \vez=\vef,
\vez\in\R^{72}_+$, as stated above. Moreover, let $G$ denote
the set of those $24$ columns of $A_{6\times 4\times 3}$ for
which $\vez^*_i>0$.

Then the set of holes in $\cone(A_{6\times 4\times 3})$ is the
set of all points that can be written uniquely as
$\sigma(\vef+\ves)$ with $\sigma\in S_6\times S_4\times S_3$
and with $\ves\in\monoid(G)$.
\end{corollary}

Finally, we have computed the Hilbert basis of the cone
associated to the semi-graphoid for $|N|=5$
\cite{Studeny-book}. It was already shown in
\cite{Hemmecke+Morton+Shiu+Sturmfels+Wienand} that the
corresponding monoid is not normal by constructing a hole via a
different method. The computation of the full Hilbert basis was
not possible at that time, neither with {\tt Normaliz}, nor
with {\tt 4ti2}. Here is the defining matrix whose columns
generate the monoid associated to the semi-graphoid for
$|N|=5$. Every $.$ corresponds to an entry $0$. $+$ and $-$
represent entries $1$ and $-1$.

{\footnotesize
\begin{verbatim}
++++++++++......................................................................
----......++++++................................................................
-...---.........++++++..........................................................
.-..-..--................++++++.................................................
..-..-.-.-............................++++++....................................
...-..-.--..............................................++++++..................
+.........---...---...+++.......................................................
.+........-..--..........---...+++..............................................
..+........-.-.-......................---...+++.................................
...+........-.--........................................---...+++...............
....+...........-..--....-..--....+++...........................................
.....+...........-.-.-................-..--....+++..............................
......+...........-.--..................................-..--....+++............
.......+..................-.-.-........-.-.-.......+++..........................
........+..................-.--..........................-.-.-.......+++........
.........+..............................-.--..............-.--............+++...
..........+.....+.....--.+.....--.--.+..........................................
...........+.....+....-.-.............+.....--.--.+.............................
............+.....+....--...............................+.....--.--.+...........
.............+............+....-.-.....+....-.-....--.+.........................
..............+............+....--.......................+....-.-....--.+.......
...............+........................+....--...........+....--.........--.+..
...................+........+.....-.-....+.....-.-.-.-.+........................
....................+........+.....--......................+.....-.-.-.-.+......
.....................+....................+.....--..........+.....--......-.-.+.
..............................+............+........--.......+........--...--..+
......................+........+..+..-......+..+..-+..--........................
.......................+........+..+.-........................+..+..-+..--......
........................+....................+..+.-............+..+.-.....+..--.
.................................+............+.....+.-.........+.....+.-..+.-.-
....................................+............+...+.-...........+...+.-..+.--
.....................................+............+...++............+...++...+++
\end{verbatim}
}

\begin{lemma}
The Hilbert basis of the cone associated to the semi-graphoid
for $|N|=5$ has $1300$ elements that come into $21$ orbits
under the underlying symmetry group $S_5\times S_2$. These are
represented by the $21$ rows of the following matrix:

{\tiny \def\h{\hspace{-0.24cm}}
\[
\begin{array}{rrrrrrrrrrrrrrrrrrrrrrrrrrrrrrrr}
0 & \h 0 & \h 0 & \h 0 & \h 0 & \h 0 & \h 0 & \h 0 & \h 0 & \h 0 & \h 0 & \h 0 & \h 0 & \h 0 & \h 0 & \h 0 & \h 0 & \h 0 & \h 0 & \h 0 & \h 0 & \h 0 & \h 0 & \h 0 & \h 0 & \h 1 & \h 0 & \h 0 & \h 0 & \h -1 & \h -1 & \h 1\\
0 & \h 0 & \h 0 & \h 0 & \h 0 & \h 0 & \h 0 & \h 0 & \h 0 & \h 0 & \h 0 & \h 0 & \h 0 & \h 0 & \h 0 & \h 1 & \h 0 & \h 0 & \h 0 & \h 0 & \h 0 & \h 0 & \h 0 & \h 0 & \h -1 & \h -1 & \h 0 & \h 0 & \h 0 & \h 0 & \h 1 & \h 0\\
2 & \h 0 & \h 0 & \h 0 & \h 0 & \h -2 & \h -1 & \h -1 & \h 1 & \h 1 & \h 0 & \h -1 & \h 1 & \h -1 & \h 1 & \h 1 & \h 1 & \h 0 & \h 0 & \h 0 & \h -1 & \h -2 & \h 1 & \h -1 & \h -1 & \h 0 & \h -1 & \h 0 & \h 1 & \h 1 & \h 0 & \h 1\\
2 & \h 0 & \h 0 & \h 0 & \h -2 & \h 0 & \h -1 & \h -1 & \h 1 & \h 0 & \h 0 & \h 1 & \h -1 & \h 1 & \h -1 & \h 1 & \h 1 & \h -1 & \h 1 & \h 0 & \h 1 & \h -1 & \h -1 & \h 1 & \h 0 & \h -1 & \h 0 & \h -2 & \h 0 & \h 0 & \h 0 & \h 2\\
1 & \h 0 & \h 0 & \h 1 & \h 0 & \h 0 & \h 0 & \h -1 & \h -1 & \h 2 & \h -1 & \h 2 & \h -1 & \h -1 & \h -1 & \h 0 & \h 0 & \h -1 & \h -1 & \h 2 & \h -1 & \h -1 & \h -1 & \h 2 & \h -1 & \h 0 & \h 0 & \h 0 & \h 1 & \h 0 & \h 0 & \h 1\\
1 & \h 0 & \h 1 & \h 1 & \h 0 & \h -1 & \h -1 & \h 1 & \h -1 & \h 1 & \h -2 & \h 0 & \h 0 & \h -1 & \h 0 & \h 1 & \h 0 & \h 1 & \h -1 & \h 0 & \h -2 & \h 0 & \h 1 & \h 1 & \h -1 & \h -1 & \h -1 & \h 1 & \h 0 & \h 1 & \h 0 & \h 1\\
0 & \h 1 & \h 1 & \h 0 & \h 1 & \h 1 & \h -1 & \h 0 & \h -2 & \h 1 & \h 0 & \h 0 & \h -2 & \h 0 & \h 0 & \h -1 & \h -1 & \h 1 & \h 0 & \h 1 & \h -2 & \h 0 & \h -1 & \h 1 & \h 1 & \h -1 & \h 0 & \h 1 & \h -1 & \h 1 & \h 0 & \h 1\\
0 & \h 1 & \h 1 & \h 0 & \h 1 & \h 1 & \h -1 & \h 0 & \h -2 & \h 0 & \h 0 & \h 0 & \h -2 & \h 0 & \h 0 & \h -1 & \h -1 & \h 1 & \h 1 & \h 1 & \h -1 & \h 1 & \h -1 & \h 1 & \h 1 & \h -1 & \h 0 & \h 0 & \h -2 & \h 0 & \h 0 & \h 2\\
1 & \h 0 & \h 0 & \h 1 & \h 0 & \h 0 & \h -1 & \h -1 & \h 0 & \h 2 & \h -1 & \h 2 & \h 0 & \h -1 & \h -1 & \h -1 & \h 1 & \h -1 & \h -1 & \h 1 & \h -1 & \h -1 & \h -1 & \h 1 & \h -1 & \h 1 & \h 0 & \h 0 & \h 1 & \h 0 & \h 0 & \h 1\\
0 & \h 1 & \h 1 & \h 1 & \h 1 & \h 1 & \h -2 & \h -1 & \h -1 & \h 1 & \h -1 & \h 1 & \h -1 & \h -1 & \h -1 & \h -2 & \h 1 & \h 0 & \h 0 & \h 1 & \h -1 & \h 0 & \h -1 & \h 1 & \h 0 & \h 1 & \h 0 & \h 0 & \h 0 & \h 0 & \h 0 & \h 1\\
0 & \h 0 & \h 1 & \h 1 & \h 1 & \h 1 & \h 0 & \h 0 & \h 0 & \h 0 & \h -1 & \h -1 & \h -1 & \h -1 & \h -1 & \h -1 & \h -1 & \h -1 & \h 1 & \h 1 & \h -1 & \h -1 & \h 1 & \h 1 & \h 1 & \h 1 & \h 0 & \h 0 & \h 0 & \h 0 & \h -2 & \h 2\\
1 & \h 1 & \h 0 & \h 1 & \h 0 & \h -1 & \h -1 & \h -1 & \h -1 & \h 0 & \h -1 & \h 1 & \h 1 & \h -1 & \h 0 & \h 1 & \h 1 & \h 0 & \h -1 & \h 1 & \h 1 & \h -1 & \h 0 & \h -1 & \h -1 & \h -1 & \h -1 & \h 0 & \h 1 & \h 0 & \h 1 & \h 1\\
0 & \h 1 & \h 0 & \h 1 & \h 1 & \h 1 & \h 0 & \h -1 & \h 1 & \h -2 & \h 0 & \h 0 & \h 0 & \h -2 & \h 1 & \h -1 & \h -1 & \h -2 & \h 1 & \h 0 & \h 0 & \h 0 & \h 1 & \h -2 & \h -1 & \h 0 & \h 1 & \h 1 & \h 1 & \h 0 & \h 1 & \h 0\\
2 & \h -1 & \h 0 & \h 0 & \h -1 & \h -1 & \h 0 & \h 0 & \h 1 & \h 1 & \h -1 & \h 1 & \h 1 & \h 1 & \h 1 & \h 0 & \h 1 & \h -2 & \h -1 & \h -1 & \h -2 & \h 0 & \h -1 & \h -1 & \h -1 & \h -1 & \h 1 & \h 1 & \h 1 & \h 1 & \h 1 & \h 0\\
0 & \h 1 & \h 1 & \h 1 & \h 1 & \h 1 & \h -2 & \h -1 & \h -1 & \h 0 & \h -1 & \h 0 & \h -1 & \h -1 & \h -1 & \h -2 & \h 1 & \h 1 & \h 1 & \h 1 & \h 0 & \h 1 & \h 0 & \h 1 & \h 1 & \h 1 & \h -1 & \h -1 & \h -2 & \h -1 & \h -1 & \h 3\\
0 & \h 1 & \h 1 & \h 1 & \h 1 & \h 1 & \h -1 & \h -1 & \h -1 & \h -1 & \h -1 & \h -1 & \h -1 & \h -1 & \h -1 & \h -1 & \h 0 & \h 0 & \h 2 & \h 2 & \h 0 & \h 0 & \h 1 & \h 1 & \h 1 & \h 1 & \h -1 & \h -1 & \h -1 & \h -1 & \h -2 & \h 3\\
1 & \h 0 & \h 0 & \h 1 & \h 2 & \h -1 & \h -1 & \h -1 & \h -1 & \h 1 & \h -1 & \h -1 & \h 1 & \h -2 & \h 0 & \h 0 & \h 2 & \h 1 & \h 0 & \h 1 & \h -1 & \h -1 & \h 1 & \h -1 & \h -1 & \h 1 & \h -2 & \h 0 & \h 0 & \h 0 & \h 0 & \h 2\\
2 & \h 0 & \h 0 & \h 0 & \h 0 & \h -2 & \h -1 & \h -1 & \h 1 & \h 2 & \h 1 & \h -1 & \h 1 & \h -1 & \h 1 & \h 1 & \h 0 & \h 0 & \h -1 & \h 0 & \h -1 & \h -2 & \h 0 & \h -2 & \h -1 & \h -1 & \h 0 & \h 1 & \h 1 & \h 1 & \h 2 & \h 0\\
3 & \h -1 & \h -1 & \h 1 & \h -1 & \h -2 & \h 0 & \h -1 & \h 0 & \h 1 & \h -1 & \h 0 & \h 1 & \h -1 & \h 1 & \h 1 & \h 1 & \h 1 & \h -1 & \h 1 & \h 0 & \h -1 & \h 1 & \h 0 & \h -1 & \h 0 & \h -2 & \h -1 & \h 1 & \h -1 & \h -1 & \h 3\\
3 & \h 1 & \h -1 & \h -1 & \h -1 & \h -2 & \h -1 & \h -1 & \h -1 & \h 2 & \h 0 & \h 0 & \h 1 & \h 0 & \h 1 & \h 1 & \h 1 & \h 1 & \h -1 & \h 1 & \h -1 & \h -1 & \h 1 & \h -1 & \h -1 & \h -1 & \h -2 & \h 0 & \h 0 & \h 0 & \h 1 & \h 2\\
2 & \h 1 & \h 0 & \h 0 & \h 0 & \h -2 & \h -1 & \h -1 & \h -1 & \h 2 & \h -1 & \h -1 & \h 1 & \h -1 & \h 1 & \h 1 & \h 1 & \h 1 & \h -1 & \h 1 & \h -1 & \h -1 & \h 2 & \h -1 & \h -1 & \h -1 & \h -2 & \h 0 & \h 0 & \h 0 & \h 1 & \h 2\\
\end{array}
\]
}
\end{lemma}

\section{Computational Approaches}

In this section we present the two computational approaches
that allowed us to solve the three challenging Hilbert basis
computations of the cones associated to $5\times 5\times
3$-tables, to $6\times 4\times 3$-tables, and to semi-graphoids
for $|N|=5$. In the first approach, we iteratively decompose
the cone into smaller cones and exploit the underlying symmetry
and set inclusion to avoid a lot of unnecessary computations.
An implementation of this approach is freely available in the
new release {\tt latte-for-tea-too-1.4} of ``LattE for tea,
too'' ({\tt http://www.latte-4ti2.de}), a joint source code
distribution of the two software packages {\tt LattE macchiato}
and {\tt 4ti2}. In the second approach, we exploit the fact
that the cones are nearly compressed; hence many cones in any
pulling triangulation are unimodular, and the same holds in
placing triangulations. Using our second approach, none of
these unimodular cones is constructed, saving a lot of
computation time. An implementation of this approach will be
freely available in the next release of Normaliz ({\tt
http://www.math.uos.de/normaliz}), together with the input
files of the examples of this paper.

\subsection{First approach: exploiting symmetry}

Let us assume that we wish to compute the Hilbert basis of a
rational polyhedral cone
$C=\cone(\ver_1,\ldots,\ver_s)\subseteq\R^n$. Moreover, assume
that $C$ has a coordinate-permuting symmetry group $S$, that
is, if $\vev\in V$ and $\sigma\in S$ then also $\sigma(\vev)\in
C$. Herein, the vector $\sigma(\vev)$ is obtained by permuting
the components of $\vev$ according to the permutation $\sigma$.

One approach to find the Hilbert basis of $C$ is to find a
regular triangulation of $C$ into simplicial cones
$C_1,\ldots,C_k$ and to compute the Hilbert bases of the
simplicial cones $C_1,\ldots,C_k$. Clearly, the union of these
Hilbert bases is a (typically non-minimal) system of generators
of the monoid of lattice points in $C$. The drawback of this
approach is that a complete triangulation of $C$ is often too
hard to accomplish.

Instead of computing a full triangulation, we compute only a
(regular) subdivision of $C$ into few cones. To this end we
remove one of the generators of the cones, say $\ver_s$,
compute the convex hull of the cone
$C'=\cone(\ver_1,\ldots,\ver_{s-1})$, and find all facets
${\mathcal F}$ of $C'$ that are visible from $\ver_s$. By
${\mathcal F}'$ we denote the set of all cones that we get as
the convex hull of a facet in ${\mathcal F}$ with the ray
generated by $\ver_s$. Then ${\mathcal F}'\cup\{C'\}$ gives a
regular subdivision of $C$, called the \emph{subdivision with
distinguished generator $\ver_s$}. Before we now subdivide
those cones in ${\mathcal F}'$ further into smaller cones, we
use the following simple observation to remove cones that can
be avoided due to the underlying symmetry given by $S$.

\begin{lemma}\label{Lemma: Dropping cones due to symmetry}
Let $C,C_1,\ldots,C_k\subseteq\R^n$ be rational polyhedral
cones such that $C=\cup_{i=1}^k C_i$ (not necessarily a
disjoint union). Suppose that there is a permutation $\sigma$
and indices $i$ and $j$ such that
$C_i\subseteq\sigma(C_j)\subseteq C$. Then the Hilbert basis of
$C$ is contained in the union of the Hilbert bases of the cones
$C_1,\ldots,C_{i-1},\sigma(C_j),C_{i+1},\ldots,C_k$.
\end{lemma}

\boproof The result follows by observing that all lattice
points in $C_i$ also belong to $\sigma(C_j)$ and thus can be
written as a nonnegative integer linear combination of the
Hilbert basis of $\sigma(C_j)$. \eoproof

If successful, this test whether $C_i$ can be dropped is a very
efficient way of removing unnecessary cones. However, the fewer
generators are present in the cones $C_1,\ldots,C_k$, the
higher the chance that this test fails. So one has to make a
trade-off between a simple test (that may fail more and more
often) and a direct treatment of each cone $C_i$. As we compute
only regular subdivisions whose cones are spanned by some of
the vectors $\ver_1,\dots,\ver_s$, each of the cones
$C_1,\ldots,C_k$ can be represented by a characteristic
$0$-$1$-vector $\chi(C_1),\ldots,\chi(C_k)$ of length $s$ that
encodes which of the generators of $C$ are present in this
cone. This makes the test $C_i\subseteq\sigma(C_j)$ comparably
cheap, as we only need to check whether
$\chi(C_i)\leq\sigma(\chi(C_j))$.

Summarizing these ideas, the symmetry exploiting approach can
be stated as follows:

\begin{enumerate}
  \item Let $C=\cone(\ver_1,\ldots,\ver_s)\subseteq\R^n$
      and ${\mathcal C}=\{C\}$.
  \item $i:=0$
  \item While ${\mathcal C}\neq\emptyset$ do
    \begin{enumerate}
      \item $i:=i+1$
      \item For all $K\in{\mathcal C}$ compute a
          subdivision with distinguished $i$th
          generator (if existent in $K$).
      \item Let ${\mathcal T}$ be the set of all cones
          in these  subdivisions.
      \item Let ${\mathcal M}$ be the set of those
          cones with a maximum number of rays.
      \item Let ${\mathcal C}\neq\emptyset$ be the set
          ${\mathcal M}$ together with all cones $T\in
          {\mathcal T}$ that are not covered by a cone
          $\sigma(M)$ with $M\in {\mathcal M}$ and
          $\sigma\in S$, see Lemma \ref{Lemma: Dropping cones due to symmetry}.
      \item Remove from ${\mathcal C}$ all simplicial
          cones and compute their Hilbert bases.
    \end{enumerate}
  \item For each computed Hilbert basis element $\veh$
      compute its full orbit $\{\sigma(\veh):\sigma\in S\}$
      and collect them in a set ${\mathcal H}$.
  \item Remove the reducible elements from ${\mathcal H}$.
  \item Return the set of irreducible elements as the
      minimal Hilbert basis of $C$.
\end{enumerate}

This quite simple approach via triangulations and elimination
of cones by symmetric covering already solves all three
presented examples. In particular, it gives a computational
proof to Lemma \ref{Lemma: 355 is normal}. The candidates for
the representatives of Hilbert basis elements can be computed
using ``LattE for tea, too'' by calling

{\footnotesize
\begin{verbatim}
dest/bin/hilbert-from-rays-symm --hilbert-from-rays="dest/bin/hilbert-from-rays"
                                      --dimension=26 S5.rays
dest/bin/hilbert-from-rays-symm --hilbert-from-rays="dest/bin/hilbert-from-rays"
                                      --dimension=43 355.short.rays
dest/bin/hilbert-from-rays-symm --hilbert-from-rays="dest/bin/hilbert-from-rays"
                                      --dimension=42 346.short.rays
\end{verbatim}
} The data files can be found on {\tt
http://www.latte-4ti2.de}. (For typographical reasons each
command has been printed in two lines.)

\subsection{Second approach: partial triangulation}

In the second approach, we build up a triangulation of the
given cone $C=\cone(\ver_1,\ldots,\allowbreak
\ver_s)\subseteq\R^n$. However, by using the following Lemma
\ref{Lemma: Removing generators} and its Corollary \ref{Lemma:
critareas}, we can avoid triangulating many regions of the
cone, since the triangulation would consist only of unimodular
cones (for which the extreme ray generators already constitute
a Hilbert basis), or, more precisely, avoid to construct
simplicial cones whose non-extreme Hilbert basis elements are
contained in previously computed simplicial cones.

In the following we describe the facets of a full-dimensional
rational cone by (uniquely determined) primitive integral
exterior normal vectors. In other words, $F=\{x\in
C:\vecc^\intercal x=0\}$ where $c$ has coprime integer entries
and $\vecc^\intercal y\le0$ for all $y\in C$.

\begin{lemma}\label{Lemma: Removing generators}
Let $C=\cone(\ver_1,\ldots,\ver_k)\subseteq\R^n$ be a rational
polyhedral cone such that
\begin{itemize}
  \item $\ver_1,\ldots,\ver_k\in\Z^n$,
  \item $\ver_1,\ldots,\ver_{k-1}$ lie in a facet of $C$
      defined by the hyperplane $\vecc^\intercal\vex=0$,
  \item $\vecc^\intercal\ver_k=1$.
\end{itemize}
Then the Hilbert basis of $C$ is the union of $\{\ver_k\}$ and
the Hilbert basis of $\cone(\ver_1,\ldots,\allowbreak
\ver_{k-1})$.
\end{lemma}

\boproof Let $\vez\in C\cap\Z^n$. Then $\vez=\sum_{i=1}^k
\lambda_i \ver_i$ for some nonnegative real numbers
$\lambda_1,\ldots,\lambda_k$. Multiplying by $\vecc^\intercal$,
we obtain
\[
\vecc^\intercal\vez=\sum_{i=1}^k \lambda_i
\vecc^\intercal\ver_i=\lambda_k\vecc^\intercal\ver_k=\lambda_k.
\]
As $\vecc,\vez\in\Z^n$, we obtain $\lambda_k\in\Z$. Hence,
$\vez$ is the sum of an nonnegative integer multiple of
$\ver_k$ and a lattice point
$\vez-\lambda_k\ver_k\in\cone(\ver_1,\ldots,\ver_{k-1})$, which
can be written as a nonnegative integer linear combination of
elements from the Hilbert basis of this cone. The result now
follows. \eoproof

This lemma implies the following fact, which excludes many
regions when searching for missing Hilbert basis elements.

\begin{corollary}\label{Lemma: critareas}
Let $\ver_1,\dots,\ver_k\in\Z^n$ such that
$C'=\cone(\ver_1,\ldots,\ver_{k-1})$ has dimension $n$, and
$C=C'+\cone(\ver_k)$. Suppose that $\ver_k\notin C'$. Moreover,
let $F_1,\dots,F_q$ be the facets of $C'$ visible from $\ver_k$
and let $\vecc_1,\ldots,\vecc_q$ the normal vectors of these
facets as introduced above. Then
\[
\Hilbert(C')\cup\{\ver_k\}\cup\bigcup\left\{\Hilbert(F_i+\cone(\ver_k)):|\vecc_i^\intercal
\ver_k|\ge 2,\ i=1,\dots,q\right\}
\]
generates $C\cap\Z^n$.
\end{corollary}

\boproof Evidently we obtain a system of generators of
$C\cap\Z^n$ if we extend the union in the proposition over all
facets $F_i$, $i=1,\dots,q$. It remains to observe that
\[
  \Hilbert(F_i+\cone(\ver_k))=\{\ver_k\}\cup \Hilbert(C'\cap F_i)
\]
if $|\vecc_i^\intercal \ver_k|=1$. But this is the statement of
Lemma \ref{Lemma: Removing generators}. \eoproof

Corollary \ref{Lemma: critareas} yields an extremely efficient
computation of Hilbert bases---provided the case
$|\vecc_i^\intercal \ver_k|\ge 2$ occurs only rarely, or, in
other words, the system $\ver_1,\dots,\ver_k$ of generators is
not too far from a Hilbert basis.

A thoroughly consequent application of Corollary \ref{Lemma:
critareas} could be realized as follows, collecting the list
$\cA(C)$ of critical simplicial cones in a recursive algorithm.
\begin{itemize}
\item[(C1)] Initially $\cA(C)$ is empty.

\item[(C2)] One searches the lexicographically first
    linearly independent subset
    $\{\ver_{i_1},\dots,\allowbreak\ver_{i_d}\}$. If the
    cone generated by these elements is not unimodular, it
    is added to $\cA(C)$.

\item[(C3)] Now the remaining elements among
    $\ver_1,\dots,\ver_s$ (if any) are inserted into the
    algorithm in ascending order. Suppose that $C'$ is the
    cone generated by the elements  processed already, and
    let $\ver_j$ be the next element to be inserted. Then
    for all facets $F_i$ of $C'$ such that
    $\vecc_i^\intercal \ver_k\ge 2$ the list $\cA(C)$ is
    augmented by $\cA(F_i+\cone(\ver_j))$.
\end{itemize}

After all the critical simplicial cones have been collected, it
remains to compute their Hilbert bases and to reduce their
union globally, together with $\{\ver_1,\dots,\ver_s\}$.

Let us add some remarks on this approach.
\begin{itemize}
\item[(a)] It is not hard to see that the list $\cA(C)$
    constitutes a subcomplex of the lexicographical
    triangulation obtained by inserting
    $\ver_1,\dots,\ver_s$. However, this fact is irrelevant
    for the computation of Hilbert bases.

\item[(b)] In an optimal list of simplicial cones each
    candidate for the Hilbert basis of $C$ would appear
    exactly once. (The candidates are the elements of the
    Hilbert bases of the simplicial cones.) The algorithm
    above cannot achieve this goal since the cones
    $F+\cone(\ver_j)$ are treated independently of each
    other. Nevertheless it yields a reasonable
    approximation.

\item[(c)] The drawback of the algorithm above is that it
    uses the Fourier-Motzkin elimination recursively for
    subcones.  Therefore Normaliz applies the algorithm
    above only on the top level and produces a full
    triangulation of the cones $F_i+\cone(\ver_k)$ for
    which $\vecc_i^\intercal \ver_k\ge 2$ (instead of the
    list $\cA(F_i+\cone(\ver_j))$).

\item[(d)] It is a crucial feature of the height $1$
    strategy that it reduces memory usage drastically.
\end{itemize}

We illustrate the size of the computation and the gain of the
improved algorithm by the data in Table \ref{data}. In the
table we use the following abbreviations: emb-dim is the
dimension of the space in which the cone (or monoid) is
embedded, dim denotes its dimension, \# rays is the number of
extreme rays, \# HB is the number of elements in the Hilbert
basis, \# full tri is the number of simplicial cones in a full
triangulation computed by \Normaliz, \# partial tri is the
number of cones in the partial triangulation, \# cand is the
number of candidates for the Hilbert basis, and \# supp hyp is
the number of support hyperplanes.

\begin{table}
\begin{center}
\begin{tabular}{|l|r|r|r|r|r|}
\hline
&&&&&semi-graph-\\
&$4\times4\times3$ & $5\times4\times3$ & $5\times5\times3$ & $6\times4\times3$ & oid $N=5$\\
 \hline
emb-dim & 40 & 47 & 55 & 54 & 32 \\
 \hline
dim  & 30 & 36 & 43 & 42 & 26\\
 \hline
\# rays & 48 & 60 & 75 & 72 & 80 \\
 \hline
\# HB & 48   &  60   &  75   &  4,392 &  1,300 \\
 \hline
\# supp hyp  & 4,948  & 29,387 & 306,955 & 153,858 & 117,978\\
 \hline
\# full tri & 2,654,000&  102,538,980 & ? & ? & ?\\
 \hline
\# partial tri & 48 & 4,320 & 775,800 & 206,064 & 3,109,495\\
 \hline
\# cand & 96 & 1,260 & 41,593 & 10,872 & 168,014\\
 \hline
\end{tabular}
\end{center}\smallskip
\caption{Data of challenging Hilbert basis
computations}\label{data}
\end{table}

In addition to the improved algorithm just presented,
parallelization has contributed substantially to the rather
short computation times that (the experimental version of)
\Normaliz\ needs for the cones considered. The computations
were done on a SUN Fire X4450 with 24 Xeon cores, but even on a
single processor machine computation times would be moderate.

\begin{remark}
Sturmfels and Sullivant \cite[3.7]{STSU} stated a very
interesting conjecture on the normality of cut monoids of
graphs without $K_5$-minors. For graphs with 7 and 8 vertices
we have used the approach via partial triangulations (and
parallelization) in order to verify the conjecture. For these
graphs no counterexample could be found.
\end{remark}

\end{document}